\documentclass[a4paper,11pt]{article}
\usepackage[latin1]{inputenc}
\usepackage{amsmath}  
\usepackage{amsfonts} 
\usepackage[english,french]{babel}
\usepackage{amssymb}
\usepackage{dsfont}\let\mathbb\mathds
\usepackage[all,dvips]{xy}
\usepackage{latexsym}
\author{Angela Ortega \\ \small{ortega@math.unice.fr}}
\title{Vari\'et\'es de Prym associ\'ees aux rev\^etements $n$-cycliques d'une courbe hyperelliptique}
\date{}
\newenvironment{Dem}{\textit{{D\'emonstration} :}}{\begin{flushright}$\Box$\end{flushright}}

\newtheorem{Pro}{Proposition}[section]
\newtheorem{Lem}{Lemme}[section]

\newcommand{\C}{JC_{\nu}}

\newcommand{\fs}{\mathcal{O}_{H}}
\newcommand{\F}{\mathop{\rm Fix}\nolimits}
\newcommand{\im}{\mathop{\rm Im}\nolimits}
\newcommand{\Nm}{\mathop{\rm N_f}\nolimits}
\newcommand{\Nmi}{\mathop{\rm N_{f_i}}\nolimits}

\newcommand{\Ker}{\mathop{\rm Ker}\nolimits}

\begin{document}

\maketitle

\section{Introduction}

Soit $H$ une courbe projective, hyperelliptique et non singuli\`ere de genre $g$ et soit $JH \simeq Pic^0(H)$
la jacobienne de $H$. On fixe un \'el\'ement $\eta \in JH$ d'ordre $n$. On note $f : C \rightarrow H$ le rev\^etement
\'etale $n$-cyclique associ\'e \`a $\eta$. Soit $JC$ la jacobienne de $C$.
\\
Soit $\Nm : JC \rightarrow JH$, $\sum n_ip_i \mapsto \sum n_if(p_i)$, l'application Norme de $f$.
On appelle vari\'et\'e de Prym $P=Prym (C/H)$ associ\'ee au rev\^etement $f$ la
composante neutre de $\Ker\Nm$. Si $\sigma : C \rightarrow C$ est l'automorphisme
d'ordre $n$ qui engendre le groupe de Galois $Gal(C/H)$, on peut \'ecrire $P= \im (1-\sigma).$
La vari\'et\'e $P$ est une sous-vari\'et\'e ab\'elienne de $JC$ de dimension $(n-1)(g-1)$ avec polarisation $\Xi$
induite par la polarisation principale de $JC$.\\
Le but de ces lignes est de montrer que lorsque $n$ n'est pas divisible par 4 la  vari\'et\'e de Prym $P$ est
isomorphe \`a un produit de jacobiennes.

\section{D\'ecomposition de la vari\'et\'e de Prym}
On pose $Y=f^*(JH)$. C'est une sous-vari\'et\'e ab\'elienne suppl\'ementaire de $P$ dans $JC$ de dimension $g$.
On peut aussi d\'ecrire $Y$ comme l'image de l'endomorphisme Norme $1+\sigma+ \cdots + \sigma^{n-1}$.
\\
\begin{Pro}
La polarisation $\Xi$ sur $P$ est du type  $(\underbrace{1,...,1}_{(n-2)(g-1)},\underbrace{n,...,n}_{g-1})$.
\end{Pro}
\begin{Dem}
Puisque $f:C \rightarrow H$ est \'etale l'application $f^*: JH \rightarrow JC$ n'est pas injective (~\cite{CAV}
11.4.3.). En fait, $\Ker f^*= \langle \eta \rangle$, un sous-groupe de $JH[n]$ d'ordre $n$. On a donc le diagramme
commutatif suivant
$$
\shorthandoff{;:!?}
\xymatrix @!0 @R=15mm @C=1.5cm{
       JH \ar[rr]^{f^*} \ar[rd]_h & & JC \\
       & JH/\langle \eta \rangle \simeq Y \ar[ru]_{i_Y} &
     }
$$
o\`u $h$ est une isog\'enie de degr\'e $n$. Soit $\Theta$ un diviseur th\^eta dans $JC$ et soit $M=\mathcal{O}_{Y}
(i_Y^* \Theta)$. On a $h^*M \simeq\mathcal{O}_{JH}(n\Theta_H) $ car $(f^*)^*\Theta \equiv n\Theta_H$
(~\cite{CAV} 12.3.1.). D'apr\`es ~\cite{Mum2} (Lemme 2, pag. 232), $K(M) \simeq \langle \eta \rangle ^{\perp}/ \langle \eta
\rangle$, o\`u $\langle \eta \rangle ^{\perp}$ est l'orthogonale de $\langle \eta \rangle$ par rapport \`a la forme
de Weil $e^{h^*M}: K(h^*M) \times K(h^*M) \rightarrow \mathbb{C}^*$. Puisque $K(h^*M)\simeq (\mathbb{Z}/n)^{2g}$ on
obtient donc $K(M) \simeq (\mathbb{Z}/n)^{2(g-1)}$. Ainsi $M$ est une polarisation du type
$(1,n, \ldots,n)$ et par ~\cite{CAV} (Cor. 12.1.5.), $\Xi = i_P^* \Theta$ est du type $(1, \dots,1,n, \ldots,n)$.
\end{Dem}
On note $i : H \rightarrow H $ l'involution hyperelliptique. Par construction on a $C = Spec ({\bf
\mathcal{A}})$, o\`u ${\bf \mathcal{A}}:=\fs \bigoplus \eta \bigoplus \cdots \bigoplus \eta^{n-1}$ est muni d'une
structure de $\fs$-alg\`ebre donn\'ee par un isomorphisme $\tau :\fs \stackrel{\sim}{\rightarrow} \eta^{n}.$ On
consid\`ere le changement de base
\begin{eqnarray}
\shorthandoff{;:!?}
\xymatrix @!0 @R=14mm @C=1.4cm{
       Spec(i^* {\bf \mathcal{A}})= & i^* C  \ar[r]^j \ar[d]_f & C \ar[d]^f &=Spec({\bf \mathcal{A}})\\
       & H \ar[r]^i & H &
     }
\end{eqnarray}
Comme l'involution $i$ agit sur $JH$ par $(-1)_{JH}$, on peut choisir un isomorphisme $\varphi: i^*\eta \stackrel{\sim}
{\rightarrow} \eta^{-1}$ de fa\c con que $\varphi^{\otimes n} = Id$ via $\tau$. On obtient ainsi un isomorphisme de
$\fs$-alg\`ebres ${\bf \mathcal{A}} \rightarrow i^*{\bf \mathcal{A}}= i^*\fs \bigoplus i^*\eta \bigoplus \cdots
\bigoplus i^*\eta^ {n-1} $. De cette fa\c con on peut identifier $i^*C$ \`a $C$ et $j$ \`a un automorphisme
v\'erifiant $j^2=1_{C}$. Observons que ce rel\`evement \`a $C$ de l'involution hyperelliptique n'est pas canonique
car il d\'epend du choix de l'isomorphisme $\varphi$.
Dans ~\cite{BL1}(Prop. 2.1.) on montre la proposition suivante

\begin{Pro}
Le rev\^etement $C \rightarrow \mathbb{P}^1$ est galoisien avec groupe de Galois $Gal(C/\mathbb{P}^1)= D_n =
\langle j,\sigma \mid j^2=\sigma^n =1, \ j\sigma j=\sigma^{-1}\rangle$.
\end{Pro}

Le groupe di\'edral $ D_n=\langle j,\sigma \rangle$ contient les involutions $j_{\nu}=j\sigma^{\nu}$ pour
$\nu=0,\ldots, n-1$. Soient $f_{\nu} :C \rightarrow  C_{\nu}:=C/\langle j_{\nu} \rangle$ les rev\^etements
doubles ramifi\'es associ\'es \`a ces involutions. Soit $g_{\nu}$ le genre de $C_{\nu}$.
On a donc le diagramme suivant

\begin{eqnarray}
\shorthandoff{;:!?}
\xymatrix{
        & C \ar[ld]_{f_0} \ar[d]^f \ar[rd]^{f_1} & \\
	C_0 \ar[rd] & H \ar[d]^{\pi} & C_1 \ar[ld]\\
        &\mathbb{P}^1 &
     }
\end{eqnarray}
Soit $W=\{x_1, \ldots, x_{2g+2} \}\subset \mathbb{P}^1$ l'ensemble des points de Weierstrass pour le rev\^etement
$\pi:H \rightarrow \mathbb{P}^1;$  on pose $S=\{ x\in W \mid \ (\pi \circ f)^{-1}(x) \textrm{ ne contient pas}$
\\ de point fixe par $j \}$ et $T=W \setminus S$.
\begin{Pro}
\hspace{1cm}\\
a) Pour $n$ impair les courbes $C_{\nu}$ sont de genre $g_{\nu}=\frac{1}{2}(n-1)(g-1)$ pour $\nu=0,\dots,n-1$.\\
b) Pour $n$ pair $ g_{\nu} = \frac{n}{2}(g-1)+1-\frac{|T|}{2}$ pour $\nu=0,\dots,n-1$.
\end{Pro}
\begin{Dem}
\hspace{1cm}\\
\emph{a)} Il suffit de montrer la proposition pour $\nu=0$. On observe que les images par $f$ des points fixes par l'
involution $j=j_0$ sont des points de ramification du rev\^etement $H \rightarrow \mathbb{P}^1$ puisque le
diagramme (1) est commutatif. Soit $q \in H$ un point fixe par $i$ et $p \in C$ un r\'elevement de $q$. Comme $f$ est
un rev\^etement non-ramifi\'e, $p$ est un point fixe par une involution $j_m$ de $C$ avec $0 \leq m \leq n$. Puisque $n$
est impair, il existe un unique $k$ modulo $n$ qui v\'erifie l'\'equation $2k \equiv m$ mod $n$. Donc, $\sigma^k p \in
f^{-1}(q)$ est le seul point fixe par $j$ dans la fibre $f^{-1}(q)$. En effet, comme $j\sigma=\sigma^{-1}j$ on a
\begin{eqnarray*}
j\sigma^k p = j\sigma^k j_m p = \sigma^{m-k} p= \sigma^{2k-k}p= \sigma^k p.
\end{eqnarray*}
En conclusion $S= \emptyset $ et on a autant de points fixes par  $j_{\nu}$ que de points de ramification du rev\^etement
$H \rightarrow \mathbb{P}^1$. Par la formule
de Hurwitz on a
\begin{eqnarray*}
g(C)-1=2(g_{\nu}-1) +g+1,
\end{eqnarray*}
et on trouve que $g_{\nu} = \frac{1}{2}(g(C)-g)= \frac{1}{2}(n(g-1)+1-g)= \frac{1}{2}(n-1)(g-1)$.\\
\emph{b)} Dans le cas $n$ pair, sur certaines fibres au-dessus des points de
ramification du rev\^etement $\pi$, l'involution $j$ n'a pas de point fixe.
On observe que si $p\in \F (j)$ alors $\sigma^{\frac{n}{2}}p \in\F (j)$
et qu'ils sont les seuls points fix\'es sur la fibre $f^{-1}(f(p))$. Donc, par la formule de Hurwitz
\begin{equation*}
g(C)-1 = 2g_{\nu}-2 + |T|
\end{equation*}
et on obtient $g_{\nu} =\frac{n}{2}(g-1) + 1- \frac{|T|}{2}$.
\end{Dem}

Les automorphismes $\sigma$ et $j$ induisent des automorphismes dans $JC$, not\'es aussi $\sigma$ et $j$. Les
applications $f_{\nu}^*: \C \rightarrow JC$  sont injectives car les rev\^etements doubles sont ramifi\'es.
On peut donc consid\'erer les jacobiennes $\C$ comme des sous-vari\'et\'es ab\'eliennes de $JC$.
Pour tout point $c \in C$ on obtient le diagramme commutatif
\begin{eqnarray}
\shorthandoff{;:!?}
\xymatrix{
	C \ar@{^{(}->}[r]^{\alpha_c} \ar[d]_{f_\nu} & JC \ar[d]^{\Nmi}  \\
	C_\nu \ar@{^{(}->}[r]^{\alpha_{f_\nu (c)}} &  JC_\nu
     }
\end{eqnarray}
pour $\nu=0,\ldots,n-1$, o\`u $\alpha_c$ est l'application de Abel-Jacobi. On peut \'ecrire $\C= \im (1+
j_{\nu})$ dans $JC$. De plus, comme $\sigma(1+j_{\nu})=(1+j_{\nu -2})\sigma$ l'automorphisme
$\sigma$ se restreint \`a des isomorphismes
\begin{eqnarray*}
\sigma: \C \rightarrow JC_{\nu-2} \qquad \textrm{pour} \quad \nu \in \mathbb{Z}/n\mathbb{Z}.
\end{eqnarray*}
On en d\'eduit que, lorsque $n$ est pair, il y a deux classes d'isomorphismes de jacobiennes qu'on note $JC_0$ et
$JC_1$; lorsque $n$ est impair toutes les jacobiennes $JC_{\nu}$ sont isomorphes.

\begin{Pro}
Les sous-vari\'et\'es $\C$ sont contenues dans la vari\'et\'e de Prym $P$.
\end{Pro}
\begin{Dem}
On a $P=\Ker (1+ \sigma +\cdots + \sigma^{n-1})^0$. D'autre part
\begin{eqnarray*}
(1+\sigma +\cdots +\sigma^{n-1})(1+j_{\nu})=  (1+\sigma +\cdots + \sigma^{n-1}+j + j_1+ \cdots +
j_{n-1})
\end{eqnarray*}
se factorise \`a travers de l'application Norme du rev\^etement $C \rightarrow \mathbb{P}^1$ et donc est l'application nulle.
En cons\'equence $\C =\im (1+j_{\nu}) \subset \Ker (1+ \sigma +\cdots +\sigma^{n-1})$, et puisque $\C$ est
connexe on a $\C \subset P$.
\end{Dem}
Soient $p_0$ et $p_1$ les projections de $JC_0 \times JC_1$ sur les deux facteurs correspondants.

\begin{Pro}
L'application $\psi = p_0+p_1 : JC_0 \times JC_1 \rightarrow P$ est un isomorphisme de vari\'et\'es
ab\'eliennes pour $n \geq 2$ non divisible par 4.
\end{Pro}
\begin{Dem}
Par la prop. 2.3 \emph{a)} dim$P=(n-1)(g-1)=$dim$(JC_0 \times JC_1)$ pour $n$ impair. Lorsque $n$ est pair on
obtient de la proposition 2.3 que dim$JC_0=g_0=\frac{n}{2}(g-1)+1-t$ et dim$JC_1 = \frac{n}{2}(g-1)+1-s$
o\`u $t+s=g+1$. On a donc dim$(JC_0 \times JC_1)=(n-1)(g-1)$. Comme $\psi$ est bien un morphisme de vari\'et\'es
ab\'eliennes il suffit de montrer que $\psi$ est injective.\\
Soit $(x,y) \in JC_0 \times JC_1$ tel que $\psi(x,y)=x+y =0$. Alors $x=(-y) \in JC_0 \cap JC_1$.
Le lemme suivant montre que n\'ecessairement $x=0$ et donc aussi $y=0$, ce qui termine la preuve. Voir~\cite{Mum}
pag. 346 pour le cas n=2.
\end{Dem}

\begin{Lem}
$ JC_0 \cap JC_1=\{ 0\}$
\end{Lem}
\begin{Dem}
\hspace{1cm}\\
\emph{a) Cas n impair.} Soit $F \in JC_0 \cap JC_1 \subset \Ker(1-j)\cap\Ker(1-j_1)\subset \F(j,\sigma)$; on a donc
$F \in \F (\sigma) \cap P \subset \F (\sigma) \cap \Ker (1+\sigma+\cdots +\sigma^{n-1}) \subset JC[n]$. Par
ailleurs, comme $\sigma^*F \simeq F$ et $f$ \'etant \'etale cela
implique que $F=f^*L$ pour un fibr\'e en droites $L \in JH$. Si de plus $j^*F \simeq F$, alors
\begin{eqnarray*}
j^*f^* L \simeq f^* i^* L \simeq f^* L.
\end{eqnarray*}
Mais dans $JH$ l'involution $i$ agit par $-1_{JH}$, i.e. $i^*L \simeq L^{-1}$. On a donc $f^*L^{-1} \simeq
(f^*L)^{-1} \simeq f^*L$ et $F=f^*L$ est un point de 2-torsion. On conclut que
$F \in \F (j,\sigma )\cap P \subset JC[2] \cap JC[n] = \{0\}$.
\\
\emph{b) Cas n=2m, m impair.} Dans ce cas on utilise le r\'esultat suivant qui est un cas particulier du Lemme
de descente d\^u \`a Kempf (~\cite{D-N}, Th\'eor\`eme 2.3):
\begin{Lem}
Soit $X$ une vari\'et\'e alg\'ebrique int\`egre sur laquelle op\`ere un groupe fini
$G$. Soit $F$ un $G$-fibr\'e vectoriel sur $X$. Alors $F$ descend \`a $X/G$ si et seulement si pour tout point
$x \in X$, le stabilisateur de $x$ dans $G$ agit trivialement sur $F_x$.
\end{Lem}
Soient $X_i:= C_i / \langle \sigma^m \rangle$ et $ X:= C/ \langle \sigma^m \rangle$, o\`u $\sigma^m$ est une
involution qui commute avec $j$ et $j_1$. On consid\`ere la tour de courbes suivante
\begin{eqnarray}
\shorthandoff{;:!?}
\def\commutatif{\ar@{}[r]|{\circlearrowleft}}
\xymatrix{
	&C \ar[ld]_{q_0} \ar[d]^q \ar[rd]^{q_1}  & \\
	C_0 \ar[d]_{r_0} \commutatif & X \ar[ld]_{f_0} \ar[d]^f \ar[rd]^{f_1} \commutatif &C_1 \ar[d]^{r_1} \\
        X_0 \ar[rd] & H \ar[d]^{\pi} &X_1 \ar[ld] \\
        & \mathbb{P}^1 &
     }
\end{eqnarray}
\\
Soit $F \in JC_0 \cap JC_1$. On sait qu'il existe des fibr\'es en droites $M_i \in JC_i, i=0,1$ tels
que $q_i^*M_i\simeq F$. Puisque $j$ et $\sigma^m$ commutent $M_0$ est invariante par $\sigma^m$. En effet,
\begin{eqnarray*}
q_0^*\sigma^{m*}M_0 \simeq \sigma^{m*}q_0^*M_0 \simeq \sigma^{m*}F \simeq F \simeq q_0^*M_0.
\end{eqnarray*}
Comme $q_0$ est ramifi\'e, $q_0^*$ est injective et donc $\sigma^{m*}M_0 \simeq M_0$. Observons que les points
de ramification du rev\^etement $r_0: C_0 \rightarrow X_0$, i.e. les points fixes par $\sigma^m$, peuvent \^etre
relev\'es aux points fixes par $j\sigma^m=j_m$ dans $C$. En effet, soit $p \in \F (\sigma^m)$ dans $C_0$ et soit
$\tilde{p} \in C$ tel que $q_0(\tilde{p})=p$. On a
\begin{eqnarray*}
q_0(\sigma^m \tilde{p})= \sigma^m q_0(\tilde{p})= \sigma^m p =p,
\end{eqnarray*}
donc $\sigma^m \tilde{p} \in q_0^{-1}(p)= \{ \tilde{p}, j\tilde{p} \}$. Comme $q: C \rightarrow X $ est non-ramifi\'e,
$\sigma^m\tilde{p} \neq \tilde{p}$. Ainsi $\sigma^m\tilde{p} = j\tilde{p}$ et on a $j\sigma^m \tilde{p}=j_m\tilde{p}=
\tilde{p}$.\\
L'action de $\sigma^m$ sur les fibres de $M_0$ au-dessus des points de
ramification de $r_0$ est la m\^eme que celle de $j_m$ sur les fibres de $F$ au-dessus des points fixes par
$j_m$ dans $C$ puisque $q_0^*M_0\simeq F$.\\
Soit $x\in \F(j_m) \subset C$, donc $x=\sigma^{\frac{m+1}{2}}y$ ou bien $x=\sigma^{\frac{3m+1}{2}}y$ o\`u $y \in
\F(j_1)$. On observe que $\langle j_m \rangle =Stab(x) $ est un sous-groupe conjug\'e de
$\langle j_1 \rangle$ qui par hypoth\`ese agit trivialement sur $F_y$, donc $j_m$ agit aussi
trivialement sur $F_x$. On en d\'eduit que $\sigma^m$ agit trivialement sur $M_{0,q_0(x)}$ et par le lemme de
descente il existe un fibr\'e $N_0 \in JX_0$ tel que $r_0^*N_0 \simeq
M_0$. De fa\c con analogue on montre l'existence d'un fibr\'e $N_1 \in JX_1$  tel que $r_1^*N_1 \simeq M_1$.\\
Comme $q^*f_0^*N_0 \simeq q^*f_1^*N_1 \simeq F$ on a
\begin{eqnarray*}
\beta:=f_0^*N_0 \otimes (f_1^*N_1)^{-1} \in \Ker q^*.
\end{eqnarray*}
Puisque $q$ est un rev\^etement double non-ramifi\'e, $\beta^2 \simeq \mathcal{O}_X$, i.e., $\beta \in JX[2]$ mais
aussi $\beta \in JX_0 \times JX_1 \simeq  Prym(X/H) $, ce dernier isomorphisme ayant \'et\'e
d\'emontr\'e dans le cas \emph{a)}.
Comme
\begin{eqnarray*}
q^*\beta \simeq \mathcal{O}_C \simeq \sigma^*\mathcal{O}_C \simeq \sigma^*q^*\beta
\simeq q^*\sigma^*\beta
\end{eqnarray*}
on a $\sigma^*\beta \in \Ker q^*$. En fait $\sigma^*\beta \simeq \beta$ et puisque $Prym(X/H) \subset \Ker (1+\sigma +\cdots
+ \sigma^{m-1})$, on a $\beta \in JX[m]\cap JX[2]=\{0\}$. Ainsi, $f_0^*N_0 \simeq f_1^*N_1 \in JX_0 \cap JX_1 =
\{ 0 \}$  par \emph{a)} et donc $F\simeq q^* \mathcal{O}_X \simeq \mathcal{O}_C$.
\end{Dem}

\section{La polarisation}

On consid\`ere la polarisation $\Xi$ dans $JC_0 \times JC_1$. On a $\Xi \equiv \psi^*\Theta$, o\`u $\Theta$ est la
polarisation principale dans $JC$. On pose $\phi =\phi_{\Xi} : JC_0 \times JC_1 \longrightarrow \widehat{JC_0}
\times \widehat{JC_1}$. Donc $\phi$ est de la forme
\begin{displaymath}
\phi = \left( \begin{array}{cc} \alpha & \hat{\beta}\\ \beta & \delta \end{array} \right)
\end{displaymath}
L'application $\alpha: JC_0 \rightarrow \widehat{JC_0}$ est la restriction \`a $JC_0$ de la polarisation principale
de $JC$. Or, l'inclusion $f_0^*: JC_0 \rightarrow JC$ est le pullback d'un rev\^etement double ramifi\'e, et
par ~\cite{CAV} (12.3.1.) on obtient $(f_0^*)^*\Theta \equiv 2\Theta_0$, o\`u $\Theta_0$ est la polarisation
principale dans $JC_0$. Ainsi $\alpha= \phi_{2\Theta_0}=2\phi_{\Theta_0}$. De fa\c con analogue, on obtient
$\delta = 2\phi_{\Theta_1}: JC_1 \rightarrow \widehat{JC_1}$, o\`u $\Theta_1$ est la polarisation principale dans
$JC_1$.\\ L'application $\beta$ est l'application qui fait commuter le diagramme
\begin{eqnarray}
\shorthandoff{;:!?}
\xymatrix{
	JC \ar[r]^{\phi_\Theta}  & \widehat{JC} \ar[d]^{\widehat{f_1^*}} \\
	JC_0 \ar@{^{(}->}[u]^{f_0^*} \ar[r]^{\beta} &  \widehat{JC_1}
     }
\end{eqnarray}
donc $\beta = \widehat{f_1^*} \circ\phi_\Theta \circ f_0^*$ et $\hat{\beta} = \widehat{f_0^*} \circ\phi_\Theta
\circ f_1^*$. Pour expliciter $\widehat{f_i^*}$ on consid\`ere le diagramme commutatif
\begin{eqnarray}
\shorthandoff{;:!?}
\xymatrix{
	JC \ar[r]^{\phi_\Theta}  & \widehat{JC} \\
	JC_i \ar@{^{(}->}[u]^{f_i^*} \ar[r]^{\phi_{\Theta_i}} &  \widehat{JC_i} \ar[u]_{\widehat{\Nmi}}
     }
\end{eqnarray}
pour $i=0,1$, obtenu en appliquant le foncteur $Pic^0$ au diagramme (3). Ensuite, en dualisant (6) on a
\begin{eqnarray}
\shorthandoff{;:!?}
\xymatrix{
	\widehat{JC} \ar[r]^{\phi_\Theta^{-1}} \ar[d]_{\widehat{f_i^*}} & JC \ar[d]^{\Nmi}\\
	\widehat{JC_i}  \ar[r]^{\phi_{\Theta_i}^{-1}} &  JC_i
     }
\end{eqnarray}
pour $i=0,1$. On a donc $\widehat{f_i^*} = \phi_{\Theta_i} \circ \Nmi \circ \phi_{\Theta}^{-1}$ et en
utilisant le fait que $\Nmi = 1+j_i$ on obtient
\begin{displaymath}
\beta = \phi_{\Theta_1} \circ (1+j_1) \circ f^*_0 \textrm{ \ \ et \ \ } \hat{\beta} = \phi_{\Theta_0}
\circ (1+j) \circ f^*_1.
\end{displaymath}

\vspace{2cm}
Angela Ortega\\
Laboratoire J.-A. Dieudonn\'e\\
Universit\'e de Nice Sophia-Antipolis \\
Parc Valrose\\
F-0610 Nice CEDEX 02, France\\
e-mail: ortega@math.unice.fr
\end{document}